\documentclass[12pt]{article}
\usepackage{amssymb}
\usepackage{amsmath}
\newcommand{\lap}{\mbox{$\bigtriangleup$}}

\newcommand{\ra}{{\mbox{$\rightarrow$}}}
\newcommand{\be}{\begin{equation}}
\newcommand{\ee}{\end{equation}}

\newtheorem{mthm}{Theorem}

\newtheorem{mcor}{Corollary}

\newtheorem{lem}{Lemma}[section]

\begin{document}

\title{\bf A semilinear equation involving the \\
fractional Laplacian in $\mathbb{R}^n$}
\medskip

\author{ Yan Li}
\maketitle

\begin{abstract}
In this paper, we consider the semilinear equation involving the fractional Laplacian in the Euclidian space $\mathbb{R}^n$:
\begin{equation}
(-\lap)^{\alpha/2} u(x) = f(x_n) \,u^p(x) , \quad x \in \mathbb{R}^n
\label{n26}
\end{equation}
in the subcritical case with $1<p<\frac{n+\alpha}{n-\alpha}$.
Instead of carrying out direct investigations on pseudo-differential equation (\ref{n26}), we first seek its equivalent form in an integral equation as below:
\begin{equation}
u(x)=\int_{\mathbb{R}^n}G_{\infty}(x,y)\,f(y_n)\, u^{p}(y)\,dy,
\label{n27}
\end{equation}
where
$ G_{\infty}(x,y)$ is the Green's function associated with the
fractional Laplacian in $\mathbb{R}^n$.
Exploiting the \emph{method of moving planes in integral forms},
we are able to derive the nonexistence of positive solutions for
(\ref{n27}) in the subcritical case. Hence the same conclusion is true for (\ref{n26}).

\end{abstract}

\bigskip


\section{Introduction}

The fractional Laplacian is a nonlocal operator
defined as
\be
(-\lap)^{\alpha/2} u(x) = C_{n,\alpha}P.V.\int_{\mathbb{R}^n}
\frac{u(x)-u(y)}{|x-y|^{n+\alpha}} dy,\label{n31}
\ee
where $\alpha$ can be any real number between 0 and 2, and P.V. stands for the Cauchy principle value.

 Besides (\ref{n31}), the fractional Laplacian has two other equivalent definitions.
One is using the {\em extension method} introduced by Caffarelli and Silvestre in \cite{CS}, and the other is by the Fourier transform:
$$\widehat{(-\lap)^{\alpha/2}}u(\xi)=|\xi|^\alpha\hat{u}(\xi),$$
with $\hat{u}$ denoting the Fourier transform of $u$.
This operator is well defined in ${\cal S}$, the Schwartz space of rapidly decreasing $C^\infty$ functions in $R^n$, and it can be extended to the distributions in the space
$${\cal L}_{\alpha/2}=\{u \mid \int_{\mathbb{R}^n}
\frac{|u(x)|}{1+|x|^{n+\alpha}}dx<\infty \}$$
by
$$<(-\lap)^{\alpha/2} u, \phi>=\int_{\mathbb{R}^n}u\,(-\lap)^{\alpha/2}\phi\, dx,
\mbox{ for all }\phi\in C_0^\infty(\mathbb{R}^n).$$

For any domain $\Omega \subset \mathbb{R}^n$ and for a given $g\in L_{loc}^1(\Omega)$, we say that $
u\in {\cal L}_{\alpha/2}$ is a solution to the problem
$$(-\lap)^{\alpha/2} u(x) =g(x), \;\; x \in \Omega$$
if and only if
\be
\int_{\mathbb{R}^n}u\,(-\lap)^{\alpha/2}\phi\, dx=
\int_{\mathbb{R}^n}g(x)\,\phi(x)\, dx,
\mbox{ for all }\phi\in C_0^\infty(\Omega).\label{n32}
\ee
In this paper, we only consider solutions in the distributional sense as given in (\ref{n32}).

We start by investigating the nonlocal equation with a specific nonlinearity in $\mathbb{R}^n$:
\begin{equation}
(-\lap)^{\alpha/2} u(x) = x_n^2\, u^p(x) , \quad x \in \mathbb{R}^n,
\label{n25}
\end{equation}
in the subcritical case with $1< p< \frac{n+\alpha}{n-\alpha}$.

Then we will
deal with the problem assuming a more general form:
\begin{equation}
(-\lap)^{\alpha/2} u(x) = f(x_n)\, u^p(x) , \quad x \in \mathbb{R}^n,
\label{n33}
\end{equation}
 where $f(x_n)$ is a positive increasing function.

 For both pseudo-differential problems (\ref{n25}) and (\ref{n33}), we want to obtain the nonexistence of
 positive solutions $u$.

 The main idea of the proof is as follows.

First, we prove that (\ref{n25}) is equivalent to the integral equation
\begin{equation}
u(x)=\int_{\mathbb{R}^n}G_{\infty}(x,y)\,y_n^2 \,u^{p}(y)\,dy,
\label{n0}
\end{equation}
where
$$ G_{\infty}(x,y) = \frac{A_{n,\alpha}}{|x-y|^{n-\alpha}}$$
is the Green's function in $\mathbb{R}^n$.

\begin{mthm}
Assume that $u \in {\cal L}_{\alpha/2} \cap L_{loc}^{\infty}(\mathbb{R}^n)$.
If $u$ is a nonnegative solution of (\ref{n25}), then $u$ satisfies the integral equation (\ref{n0}), and vice versa.
\label{mthm1}
\end{mthm}

Thanks to the equivalence between (\ref{n25}) and (\ref{n0}), in order to verify that (\ref{n25}) admits no positive solutions, it suffices
to show that same conclusion holds for the integral equation (\ref{n0}). Utilizing the \emph{method of moving planes in integral forms}, we verify

\begin{mthm}
   Assume that $u \in {\cal L}_{\alpha/2} \cap L_{loc}^{\infty}(\mathbb{R}^n)$.
   If $u$ is a nonnegative  solution of (\ref{n0}) when $1< p< \frac{n+\alpha}{n-\alpha}$,
 then $u(x)\equiv 0$.
 \label{mthm2}
\end{mthm}

As an immediate consequence, we have the following
\begin{mcor}
 Assume that $u\in {\cal L}_{\alpha/2} \cap
  L_{loc}^{\infty}(\mathbb{R}^n)$.
   If $u$ is a nonnegative solution of (\ref{n25}) when $1< p< \frac{n+\alpha}{n-\alpha}$,
 then $u(x)\equiv 0$.
\end{mcor}

For (\ref{n33}), similarly, we obtain the equivalence in
the first place.

\begin{mthm}
Assume that $u \in {\cal L}_{\alpha/2} \cap L_{loc}^{\infty}(\mathbb{R}^n)$
 and $f$ is an increasing positive function.
If $u$ is a nonnegative solution of (\ref{n33}), then $u$ satisfies the integral equation
\be
u(x)=\int_{\mathbb{R}^n}G_{\infty}(x,y)\,f(y_n)\,u^{p}(y)\,dy,
\label{n34}
\ee
and vice versa.
\label{mthm3}
\end{mthm}

Then we employ the \emph{method of moving planes} to show that
\begin{mthm}
   Assume that $u\in {\cal L}_{\alpha/2} \cap L_{loc}^{\infty}(\mathbb{R}^n)$ and $f$ is an increasing positive function. If $u$ is a nonnegative solution of (\ref{n34}) when $1< p< \frac{n+\alpha}{n-\alpha}$,
 then $u(x)\equiv 0$.
 \label{mthm4}
\end{mthm}

Instantly, we arrive at

\begin{mcor}
 Assume that $u\in {\cal L}_{\alpha/2} \cap
  L_{loc}^{\infty}(\mathbb{R}^n)$ and $f$ is an increasing positive function. If $u$ is a nonnegative solution of (\ref{n33}) when $1< p< \frac{n+\alpha}{n-\alpha}$,
 then $u(x)\equiv 0$.
\end{mcor}

The significance of such results on the nonexistence of global positive solutions lies in the fact that it serves as an important ingredient in obtaining a prior estimate for solutions of a corresponding family of nonlocal equations on bounded domains in the Euclidean space or on Riemannian manifolds. For more articles concerning the Liouville type theorem and integral
equations, please see \cite{CDL}, \cite{CFY}, \cite{CL2}, \cite{LZ},
\cite{RFB}, \cite{ZW}, \cite{ZCCY} and the references therein.

The paper is organised as follows: we present the main results in the introduction. The second section is devoted to the proof of Theorem
\ref{mthm1} and Theorem \ref{mthm2}. In the third section, we
verify Theorem \ref{mthm3} and Theorem \ref{mthm4}.

\section{Nonexistence of positive solution for \\
$(-\lap)^{\alpha/2} u(x) = x_n^2\, u^p(x) $}

\subsection{Equivalence between the integral equation and the pseudo-differential equation}

In \cite{Ku}, the Green's function of the operator $(-\lap)^{\alpha/2}$ with the Dirichlet condition on the unit ball $B_1$ is obtained as:
$$
G_{1}(x,y)= \frac{A_{n,\alpha}}{s^{\frac{n-\alpha}{2}}}\left[1-
B_{n,\alpha}\frac{1}{(t+s)^{\frac{n-2}{2}}}
\int_{0}^{\frac{s}{t}}\frac{(s-tb)^{\frac{n-2}{2}}}{b^{\alpha/2}(1+b)}db\right]
, x,y \in B_1,
$$
where $s=|x-y|^{2}$, $t=(1-|x|^{2})(1-|y|^{2})$ and $A_{n,\alpha}$
and $B_{n,\alpha}$ are constants relying on $n$ and $\alpha$.

Then we
can write the Green's function on $B_{R}$ as
\begin{eqnarray}
\nonumber
&&G_{R}(x,y)\\
\nonumber
&=&\frac{1}{R^{n-\alpha}}G_{1}(\frac{x}{R},\frac{y}{R})\\
\label{n1}
&=&\frac{A_{n,\alpha}}{|x-y|^{n-\alpha}} \left[1-
B_{n,\alpha} \frac{1}{\left(1+\frac{s_{R}}{t_{R}}\right)^{\frac{n-2}{2}}}
\int_{0}^{\frac{s_{R}}{t_{R}}}\frac{\left(\frac{s_{R}}{t_{R}}-b
\right)^{\frac{n-2}{2}}}{b^{\alpha/2}(1+b)}db\right],
\end{eqnarray}
with $s_{R}=\frac{|x-y|^{2}}{R^{2}}$ and $t_{R}=\left(1-|\frac{x}{R}|^{2}\right)\left(1-|\frac{y}{R}|^{2}\right)$.

Let $$v_R(x)=\int_{B_R}G_{R}(x,y)\,y_n^2 \,u^{p}(y)\,dy.$$

Since $u\in L_{loc}^{\infty}(\mathbb{R}^n)$, it is obvious that for each $R>0$,
$v_R(x)$ is well-defined and satisfies
\begin{equation}
\left\{\begin{array}{ll}
(-\lap)^{\alpha/2} v_R(x) = x_n^2 \,u^p(x) , & x \in B_R , \\
v_R(x) = 0 , & x \not\in B_R .
\end{array}
\right.
\end{equation}

Let $w_R(x)=u(x)-v_R(x)$ and we have
\begin{equation}
\left\{\begin{array}{ll}
(-\lap)^{\alpha/2} w_R(x) = 0 , & x \in B_R , \\
w_R(x)\geq 0 , & x \not\in B_R .
\end{array}
\label{n2}
\right.
\end{equation}

To continue, we need the \emph{maximum principle
} for factional
Laplacians.
\begin{lem}(\cite{Si})
Let $\Omega \subseteq \mathbb{R}^n$ be an open bounded set, and let $u$ be
a lower-semicontinuous function in $\overline{\Omega}$ such that $(-\lap)
^{\alpha/2}u\geq 0$ in $\Omega$ and $u\geq0$ in $\mathbb{R}^n\backslash \Omega$.
Then $u\geq 0$ in $\mathbb{R}^n$.
\label{mp}
\end{lem}
Applying Lemma \ref{mp} to (\ref{n2}), we arrive at
$$w_R(x)\geq 0, \quad x \in B_R.$$
Thus
$$w_R(x)\geq 0 , \quad x \in \mathbb{R}^n.$$

Sending $R\rightarrow\infty$ in (\ref{n1}) and we obtain the Green's function $G_{\infty}(x,y)$ in $\mathbb{R}^n$:
$$ G_{\infty}(x,y) = \frac{A_{n,\alpha}}{|x-y|^{n-\alpha}}.$$
Meanwhile,
$$v_R(x) \ra v(x):=\int_{\mathbb{R}^n}\frac{A_{n,\alpha}}{|x-y|^{n-\alpha}}\,y_n^2 \,u^{p}(y)\,dy,$$
and
$$w_R(x) \ra w(x):=u(x)-v(x).$$
It follows from (\ref{n2}) that
\begin{equation}
\left\{\begin{array}{ll}
(-\lap)^{\alpha/2} w(x) = 0 , & x \in \mathbb{R}^n , \\
w(x)\geq 0 , & x \in \mathbb{R}^n.
\end{array}
\label{n3}
\right.
\end{equation}

\begin{lem}(\cite{ZCCY})
Every $\alpha$-harmonic function bounded either above or below in all
of $\mathbb{R}^n$ for $n\geq2$ must be constant.
\label{lt}
\end{lem}
The above \emph{Liouville theorem} for fractional Laplacians implies that
$$w(x) \equiv C \geq 0 , \quad x \in \mathbb{R}^n.$$
 Here and below $C$ stands for nonnegative constants of different values 
 in various line.
 Now it's obvious that $v(x)$ is well-defined. To establish the equivalence, we need to show that $C=0$. Indeed, if $C>0$, then for each fixed $x\in \mathbb{R}^n$,
\begin{eqnarray*}
\infty>u(x)&\geq & v(x)\\
&=&\int_{\mathbb{R}^n}\frac{A_{n,\alpha}}{|x-y|^{n-\alpha}}\,y_n^2\, u^{p}(y)\,dy\\
&\geq& \int_{\mathbb{R}^n}\frac{A_{n,\alpha}}{|x-y|^{n-\alpha}}\,y_n^2 \,C^p\,dy\\
&\geq& C\int_{\mathbb{R}^n\backslash D}\frac{ y_n^2}{|x-y|^{n-\alpha}}dy\\
&\geq& C\int_{\mathbb{R}^n\backslash D}\frac{dy}{|x-y|^{n-\alpha}}\\
&=&\infty,
\end{eqnarray*}
with $D=\{y\in \mathbb{R}^n\mid |y_n|<1\}$.
The contradiction above implies that $w=0$. Therefore,
$$u(x)=v(x)=\int_{\mathbb{R}^n}\frac{A_{n,\alpha}}{|x-y|^{n-\alpha}}\,y_n^2 \,u^{p}(y)\,dy.$$

Next we prove that if $u(x)$ solves the integral equation (\ref{n0}), it also solves the differential equation. For any $\phi \in C_0^\infty(\mathbb{R}^n)$,
\begin{eqnarray*}
<(-\lap)^{\frac{\alpha}{2}}u, \, \phi>&=&<\int_{\mathbb{R}^n}G_{\infty}(x,y)\,y_n^2 \,u^{p}(y)\,dy, \, (-\lap)^{\frac{\alpha}{2}}\phi(x)>\\
&=& \int_{\mathbb{R}^n}\left\{\int_{\mathbb{R}^n}G_{\infty}(x,y)\,y_n^2\, u^{p}(y)\,dy\right\}(-\lap)^{\frac{\alpha}{2}}\phi(x)\,dx\\
&=& \int_{\mathbb{R}^n}\left\{\int_{\mathbb{R}^n}G_{\infty}(x,y) (-\lap)^{\frac{\alpha}{2}}\phi(x)dx\right\}\,y_n^2\,u^{p}(y)\,dy\\
&=& \int_{\mathbb{R}^n}\left\{\int_{\mathbb{R}^n}\delta(x-y) \phi(x)dx\right\}\,y_n^2\,u^{p}(y)\,dy\\
&=& \int_{\mathbb{R}^n}\phi(y)\,y_n^2\,u^{p}(y)\,dy\\
&=& <y_n^2\,u^{p},\, \phi>.
\end{eqnarray*}
Hence $u(x)$ satisfies (\ref{n25}). This proves Theorem \ref{mthm1}.

\subsection{Nonexistence of positive solutions for the integral equation}

\quad We will use the Kelvin type transform and the method of moving planes to derive the nonexistence of positive solutions of (\ref{n0}) under the assumption that $u\in {\cal L}_{\alpha/2}\cap
L_{loc}^{\infty}(\mathbb{R}^n)$.

With no global integrability assumptions on the solution of (\ref{n0}), we cannot carry out even the first step in the method of moving planes. To overcome this difficulty, we turn to $\bar{u}$, the Kelvin type transform of $u$, which has the desired integrability at infinity.

For $z^0=({z^0}^\prime, 0)\in (\mathbb{R}^{n-1}, \mathbb{R})$, let
\be
\bar{u}(x)=\frac{1}{|x-z^{0}|^{n-\alpha}}
u\left(\frac{x-z^{0}}{|x-z^{0}|^{2}}+z^{0}\right)
\label{n4}
\ee
be the Kelvin type transform of $u$ centered at $z^0$. Apparently,
$\bar{u}$ is integrable near infinity.

Through an elementary calculation we have
\begin{eqnarray}
\nonumber
\bar{u}(x)
&=&\frac{1}{|x-z^{0}|^{n-\alpha}}
\int_{\mathbb{R}^n}G_{\infty}\left(\frac{x-z^{0}}
{|x-z^{0}|^{2}}+z^{0},y \right)\,y_n^2\, u^{p}(y)\,dy \nonumber \\
&=&\int_{\mathbb{R}^n}G_{\infty}(x,y)\frac{y_n^2 \, \bar{u}^{p}(y)}{|y-z^{0}|^{\beta}}dy,
\label{n5}
\end{eqnarray}
for all $x\in \mathbb{R}^n\backslash B_{\epsilon}(z^{0})$ and
$\epsilon>0$, where $\beta=4+(n-\alpha)(\tau-p)\geq 0$ and $\tau=\frac{n+\alpha}{n-\alpha}$.

Now we carry out the method of moving planes on a nonnegative solution $\bar{u}$ of (\ref{n5}). Our goal is
to show that $\bar{u}$ is symmetric about the line passing through $z^0$ and parallel to the $x_n$ axis. Such symmetry enables us to
derive that $u$ is independent of its first $(n-1)$th variables
$x_1$, $\cdots$, $x_{n-1}$, and consequently obtain that
$$u=u(x_n).$$
The fact that the value of $u(x)$ is determined by its $x_n$th variable only will lead to a contradiction with the finiteness of the integral
$$\int_{\mathbb{R}^n}G_{\infty}(x,y)\,y_n^2 \,u^{p}(y)\,dy.$$ By then it's easy
to see that $u(x)$ must be trivial.

We begin the proof by introducing some notations. For a given real number $\lambda$, let the moving plane be
$$T_\lambda=\{x \in \mathbb{R}^n\mid x_1=\lambda\}.$$
Let
$$
\Sigma_\lambda=\{x=(x_{1},\cdots,x_{n})\in
\mathbb{R}^n \mid x_{1}<\lambda\},
$$
and
$$
x^{\lambda}=(2\lambda-x_{1},x_{2},\cdots,x_{n})
$$
be the reflection of the point $x$ about $T_\lambda$. Set
$$u_\lambda(x)=u(x^\lambda) \mbox{ and }w_\lambda(x)=\bar{u}_\lambda(x)-\bar{u}(x).$$

The argument will be presented in two parts. In the first part,
we begin moving the plane
$T_\lambda$ from the neighbourhood of $x_1= -\infty$. We want to show that for $\lambda$ sufficiently negative,
\be
w_\lambda(x)\geq0, \mbox{ a.e. in } \Sigma_\lambda.
\label{n28}
\ee
This provides a starting point for the moving of the planes. As long as (\ref{n28}) holds, we can keep moving the planes to the right until it reaches a limiting position $\lambda=z_1^0$. Going through a similar argument, one can move $T_\lambda$ from the positive infinity to the left and show that
\be
w_\lambda(x)\leq0, \mbox{ a.e. in } \Sigma_\lambda,
\label{n29}
\ee
for any $\lambda\geq z_1^0$. Combining (\ref{n28}) and (\ref{n29}), it's trivial to obtain the symmetry of $u_\lambda(x)$
about the plane $T_{z_1^0}$, i.e.
\be
w_{z_1^0}\equiv0,  \mbox{ a.e. in } \Sigma_{z_1^0}.
\ee
This achieves the goal in part two.

Here is the detailed proof.

\emph{Step 1.} \emph{Start moving the planes from $-\infty$ to the right as long as (\ref{n28}) holds.}

For any $\epsilon>0$, define
$$
\Sigma_{\lambda}^{-}=\{x\in \Sigma_{\lambda}\backslash B_{\epsilon}((z^{0})^{\lambda})\mid w_{\lambda}(x)<0\}.
$$
We show that for $\lambda$ sufficiently negative,
$\Sigma_{\lambda}^{-}$ must be measure zero.

By (\ref{n5}), we have
$$
\bar{u}(x)=\int_{\Sigma_{\lambda}}G_{\infty}(x,y)\frac{y_n^2\, \bar{u}^{p}(y)}{|y-z^{0}|^{\beta}}dy
+\int_{\Sigma_{\lambda}}G_{\infty}(x,y^{\lambda})
\frac{y_n^2\,\bar{u}^{p}_{\lambda}(y)}{|y^{\lambda}-z^{0}|^{\beta}}dy,
$$
$$
\bar{u}(x^{\lambda})
=\int_{\Sigma_{\lambda}}G_{\infty}(x^{\lambda},y)
\frac{y_n^2\,\bar{u}^{p}(y)}{|y-z^{0}|^{\beta}}dy
+\int_{\Sigma_{\lambda}}G_{\infty}(x^{\lambda},y^{\lambda})
\frac{y_n^2\,\bar{u}^{p}_{\lambda}(y)}{|y^{\lambda}-z^{0}|^{\beta}}dy.
$$
Then
\begin{eqnarray}
&&\bar{u}(x)-\bar{u}_{\lambda}(x)\nonumber\\
&=&\int_{\Sigma_{\lambda}}\left[G_{\infty}(x,y)-G_{\infty}(x^{\lambda},y)\right]
\frac{y_n^2 \,\bar{u}^{p}(y)}{|y-z^{0}|^{\beta}}dy \nonumber\\ &+&\int_{\Sigma_{\lambda}}
\left[G_{\infty}(x,y^{\lambda})-G_{\infty}
(x^{\lambda},y^{\lambda})\right]\frac{y_n^2 \, \bar{u}^{p}_{\lambda}(y)}{|y^{\lambda}-z^{0}|^{\beta}}dy\nonumber\\
&=&\int_{\Sigma_{\lambda}}\left[
G_{\infty}(x,y)-G_{\infty}
(x^{\lambda},y)\right]\left[\frac{y_n^2\,\bar{u}^{p}(y)}{|y-z^{0}|
^{\beta}}-\frac{y_n^2\,\bar{u}_{\lambda}^{p}(y)}{|y^{\lambda}-z^{0}|^{\beta}}\right]dy.
\label{n5.1}
\end{eqnarray}

By the Mean Value Theorem, for sufficiently negative values of $\lambda$ and $x\in\Sigma_{\lambda}^{-}$, we have
\begin{eqnarray}
\nonumber
0&<&\bar{u}(x)-\bar{u}_{\lambda}(x)\\
\nonumber
&= &\int_{\Sigma_{\lambda}}y_n^2\left[
G_{\infty}(x,y)-G_{\infty}
(x^{\lambda},y)\right]\left[\frac{\bar{u}^{p}(y)}{|y-z^{0}|^{\beta}}
-\frac{\bar{u}_{\lambda}^{p}(y)}{|y^{\lambda}-z^{0}|^{\beta}}\right]dy\\
\nonumber
&=&\int_{\Sigma_{\lambda}^{-}}y_n^2\left[
G_{\infty}(x,y)-G_{\infty}
(x^{\lambda},y)\right]\left[\frac{\bar{u}^{p}(y)}{|y-z^{0}|^{\beta}}
-\frac{\bar{u}_{\lambda}^{p}(y)}{|y^{\lambda}-z^{0}|^{\beta}}\right]dy\\
\nonumber
&\;\;+&\int_{\Sigma_{\lambda}\backslash \Sigma_{\lambda}^{-}}y_n^2\left[
G_{\infty}(x,y)-G_{\infty}
(x^{\lambda},y)\right]\left[\frac{\bar{u}^{p}(y)}{|y-z^{0}|^{\beta}}
-\frac{\bar{u}_{\lambda}^{p}(y)}{|y^{\lambda}-z^{0}|^{\beta}}\right]dy
\\
\nonumber
&\leq& \int_{\Sigma_{\lambda}^{-}}y_n^2\left[
G_{\infty}(x,y)-G_{\infty}
(x^{\lambda},y)\right]\left[\frac{\bar{u}^{p}(y)}{|y-z^{0}|^{\beta}}
-\frac{\bar{u}_{\lambda}^{p}(y)}{|y^{\lambda}-z^{0}|^{\beta}}\right]dy
\\
\nonumber
&=&\int_{\Sigma^{-}_{\lambda}}y_n^2\left[
G_{\infty}(x,y)-G_{\infty}(x^{\lambda},y)\right]
\bigg[\frac{\bar{u}^{p}(y)-\bar{u}_{\lambda}^{p}(y)}
{|y-z^{0}|^{\beta}}\\\nonumber
&&+\bar{u}_{\lambda}^{p}(y)
[\frac{1}{|y-z^{0}|^{\beta}}-\frac{1}{|y^{\lambda}-z^{0}|^{\beta}}]\bigg]dy\\
\nonumber
&\leq&\int_{\Sigma^{-}_{\lambda}}y_n^2\left[
G_{\infty}(x,y)-G_{\infty}(x^{\lambda},y)\right]
\frac{\bar{u}^{p}(y)-\bar{u}_{\lambda}^{p}(y)}
{|y-z^{0}|^{\beta}}dy\\
\nonumber
&\leq&p\int_{\Sigma^{-}_{\lambda}}y_n^2
G_{\infty}(x,y)\frac{\bar{u}^{p-1}(y)}{|y-z^{0}|^{\beta}}
[\bar{u}(y)-\bar{u}_{\lambda}(y)]dy\\
\label{n6}
&\leq&\int_{\Sigma^{-}_{\lambda}}
\frac{C}{|x-y|^{n-\alpha}}\,\frac{\bar{u}^{p-1}(y)}{|y-z^{0}|^{\beta-2}}\,
(\bar{u}(y)-\bar{u}_{\lambda}(y))dy.
\end{eqnarray}

To continue, we need the following lemma.
\begin{lem}
\label{le1} (An equivalent form of the Hardy-Littlewood-Sobolev inequality)  Assume $0<\alpha<n$ and $\Omega\subset \mathbb{R}^n$. Let $g\in L^{\frac{np}{n+\alpha p}}(\Omega)$ for $\frac{n}{n-\alpha}<p<\infty.$ Define
$$
Tg(x):=\int_{\Omega}\frac{1}{|x-y|^{n-\alpha}}g(y)dy.
$$
Then
\begin{equation}\label{n8}
 \|Tg\|_{L^{p}(\Omega)}\leq C(n,p, \alpha)\|g\|_{L^{\frac{np}{n+\alpha p}}(\Omega)}.
\end{equation}
 \end{lem}
The proof of this lemma is standard and can be found in \cite{CL1} or \cite{CL2}.

Let $\Omega$ be any domain that is a positive distance away from $z^0$. Since $u$ is locally bounded, we have
\begin{equation}\label{n7}
\int_{\Omega}\left[\frac{\bar{u}^{p-1}(y)}{|y-z^{0}|^{\beta-2}}\right]
^{\frac{n}{\alpha}}dy<\infty.
\end{equation}

For any $q>\frac{n}{n-\alpha}$, applying the Hardy-Littlewood-Sobolev inequality (\ref{n8}) and H\"{o}lder
inequality to (\ref{n6}) yields
\begin{eqnarray}\nonumber
\|w_{\lambda}\|_{L^{q}(\Sigma_{\lambda}^{-})}&\leq&
C\|\frac{\bar{u}^{p-1}}{|y-z^{0}|^{\beta-2}}w_{\lambda}\|_{L^{\frac{nq}{n+\alpha
q}}(\Sigma_{\lambda}^{-})}\\
\label{n9}
&\leq& C\|\frac{\bar{u}^{p-1}}{|y-z^{0}|^{\beta-2}}\|_{L^{\frac{n}{\alpha}}
(\Sigma_{\lambda}^{-})}
\|w_{\lambda}\|_{L^{q}(\Sigma_{\lambda}^{-})}.
\end{eqnarray}
When $N$ is sufficiently large, (\ref{n7}) indicates that for $\lambda\leq-N$,
\begin{equation}\label{n10}
C\left
\{\int_{\Sigma_{\lambda}^{-}}
[\frac{\bar{u}^{p-1}}{|y-z^{0}|^{\beta-2}}]^{\frac{n}{\alpha}}dy\right\}
^{\frac{\alpha}{n}}\leq\frac{1}{2}.
\end{equation}
A combination of (\ref{n9}) and (\ref{n10}) gives
$$
\|w_{\lambda}\|_{L^{q}(\Sigma_{\lambda}^{-})}=0.
$$
Therefore $\Sigma_{\lambda}^{-}$ must be measure zero, that is,
\begin{equation}\label{n11}
w_{\lambda}(x)\geq 0, \;\;a.e.\;x\in \Sigma_{\lambda}.
\end{equation}

\emph{Step 2.} {\em Move the plane to the limiting position to
derive symmetry}.

Inequality (\ref{n11}) serves as a starting point to move the planes $T_{\lambda}$. As long as (\ref{n11}) is valid, we will continue moving the planes to the right until the limiting position. Define
$$
\lambda_{0}=\sup\{\lambda\leq z_{1}^{0}\mid w_{\rho}(x)\geq 0,\; \rho\leq\lambda,\;
\forall x\in\Sigma_{\rho}\}.
$$

We prove
\be
\lambda_{0}\geq z_{1}^{0} - \epsilon
\label{n11.1}
\ee
using the contradiction argument.
If (\ref{n11.1}) is not true, or, $\lambda_{0}<z_{1}^{0} -\epsilon$, then
we are able to show that $\bar{u}(x)$ is symmetric about the plane
$T_{\lambda_{0}}$, i.e.
\begin{equation}\label{n12}
w_{\lambda_{0}}\equiv 0,\;\;\mbox{ a.e. in }
\Sigma_{\lambda_{0}}.
\end{equation}
If (\ref{n12}) does not hold, then  for such  $\lambda_{0} <z_{1}^{0} -\epsilon$,
$$w_{\lambda_{0}}> 0 ,
\mbox{ a.e. in } \Sigma_{\lambda_{0}}.$$
 This enables us to move the plane even further to the right. More rigorously, there exists a
$\zeta>0$ such that for all $\lambda\in[\lambda_{0},\lambda_{0}+\zeta)$,
$$w_{\lambda}(x) \geq 0 , \mbox{ a.e. in }\Sigma_{\lambda}.
$$
This contradicts the definition of $\lambda_0$.

From inequality (\ref{n9}) it follows that
\begin{equation}\label{n30}
\|w_{\lambda}\|_{L^{q}(\Sigma_{\lambda}^{-})}\leq
C\left\{\int_{\Sigma^{-}_{\lambda}}
 [\frac{\bar{u}^{p-1}(y)}{|y-z^{0}|^{\beta-2}}]^{\frac{n}{\alpha}} dy\right\}^{\frac{\alpha}{n}}
\|w_{\lambda}\|_{L^{q}(\Sigma_{\lambda}^{-})}.
\end{equation}
 When $\zeta$ is sufficiently small, for all $\lambda\in[\lambda_{0},\lambda_{0}+\zeta)$,
\begin{equation}\label{n13}
C \left\{\int_{\Sigma^{-}_{\lambda}}
[\frac{\bar{u}^{p-1}(y)}{|y-z^{0}|^{\beta-2}}]^
{\frac{n}{\alpha}}dy \right\}^{\frac{\alpha}{n}}\leq\frac{1}{2}.
\end{equation}
We will give the proof of the inequality above later. For now, from
(\ref{n30}) and (\ref{n13}), we can deduce that
$\|w_{\lambda}\|_{L^{q}_{\Sigma_{\lambda}^{-}}}=0$. Hence
$\Sigma_{\lambda}^{-}$ must be measure zero. And for
 $\lambda>\lambda_{0}$, we have
$$
w_{\lambda}(x)\geq 0,\;\;\mbox{ a.e. in }
\Sigma_{\lambda}.
$$
This is a contradiction with the definition of $\lambda_{0}$. Therefore
(\ref{n12}) holds. So far, we have verified that if $\lambda_{0}<z_{1}^{0} -\epsilon$ for any $\epsilon>0$, then
$$
\bar{u}(x)\equiv \bar{u}_{\lambda_{0}}(x),\;\;\mbox{ a.e. in }
\Sigma_{\lambda_{0}}.
$$
Furthermore, due to the singularity of $\bar{u}$ at $z^{0}$, $\bar{u}$ is also singular at $(z^{0})^{\lambda_{0}}$. It cannot be true given that $z^0$ is the only singularity of $\bar{u}$ from (\ref{n4}). This proves that
\be
\lambda_0 \geq z_1^0 - \epsilon.
  \ee
Together with the arbitrariness of $\epsilon>0$, it implies that
$$
w_{z_1^{0}}(x)\geq 0,\;\;\mbox{ a.e. in }\Sigma_{z_1^{0}}.
$$
Similarly, we can move the plane from near $x_{1}=+\infty$ to the left and show that
\be
w_{z_1^{0}}(x)\leq 0, \;\;\mbox{ a.e. in }\Sigma_{z_1^{0}}.
\ee
 Therefore,
\be\label{n13.1}
w_{z_1^{0}}(x)\equiv 0,\;\;\;\;\mbox{ a.e. in }\Sigma_{z_1^{0}}.
\ee

Now we prove inequality (\ref{n13}). For any positive small $\eta$ and
$\epsilon$, for $R$ sufficiently large it holds
\begin{equation}\label{n14}
\left(\int_{(\mathbb{R}^n\backslash B_{\epsilon}((z^{0})^\lambda))\backslash
B_{R}} \left(\frac{\bar{u}^{p-1}(y)}{|y-z^{0}|^{\beta-2}}\right)
^{\frac{n}{\alpha}}
dy\right)^{\frac{\alpha}{n}}<\eta.
\end{equation}
For any fixed $R$ large we show that the measure of
$\Sigma_{\lambda}^{-}\cap B_{R}$ is sufficiently small for
$\lambda$ close to $\lambda_{0}$. Actually, we have
\begin{equation}\label{n15}
w_{\lambda_{0}}(x)>0
\end{equation}
in the interior of $\Sigma
_{\lambda_{0}}\backslash B_{\epsilon}((z^{0})^{\lambda_{0}})$.

Indeed, if (\ref{n15}) is violated, then there exists some point
$x_0 \in \Sigma_{\lambda_0}$ such that $u(x_0)=u_{\lambda_0}(x_0)$.
From (\ref{n5.1}) we have
\begin{eqnarray*}
0&=&\bar{u}(x_0)-\bar{u}_{\lambda_0}(x_0)\\
&=&\int_{\Sigma_{\lambda_0}}y_n^2\left[
G_{\infty}(x_0,y)-G_{\infty}
(x_0^{\lambda_0},y)\right]\left[\frac{\bar{u}^{p}(y)}{|y-z^{0}|
^{\beta}}-\frac{\bar{u}_{\lambda_0}^{p}(y)}{|y^{\lambda_0}-z^{0}|^{\beta}}\right]dy.
\end{eqnarray*}
And further
$$u(x) > u_{\lambda_0}(x), \quad\forall x \in \Sigma_{\lambda_0}.$$
This obviously contradicts the fact that $w_{\lambda_0}(x)\geq 0$ in
$\Sigma_{\lambda_0}$. Hence (\ref{n15}) is true.

On the other hand, the measure of
$(\Sigma^{-}_{\lambda}\backslash B_\epsilon((z^0)^\lambda))\cap B_{R}$ can also be made as small as we want. Combining this with (\ref{n14}), (\ref{n13}) follows easily.

For any $\gamma>0$, let
$$E_{\gamma}=\{(x\in \Sigma_{\lambda_{0}}\backslash B_\epsilon((z^0)^{\lambda_{0}}))\cap
B_{R}\mid w_{\lambda_{0}}(x)>\gamma\},$$
$$F_{\gamma}=((\Sigma_{\lambda_{0}}\backslash B_\epsilon((z^0)^{\lambda_{0}}))\cap
B_{R})\backslash E_{\gamma}.
$$
Obviously,
$$
\lim_{\gamma\rightarrow 0}\mu(F_{\gamma})=0.
$$
For $\lambda>\lambda_{0}$, let
$$
D_{\lambda}=(((\Sigma_{\lambda}\backslash B_\epsilon((z^0)^\lambda))\backslash
(\Sigma_{\lambda_{0}}\backslash B_\epsilon((z^0)^{\lambda_{0}}))))\cap
B_{R}.
$$
Then
\begin{equation}
\label{n16}
((\Sigma^{-}_{\lambda}\backslash B_\epsilon((z^0)^\lambda))\cap B_{R})\subset ((\Sigma^{-}_{\lambda}\backslash B_\epsilon((z^0)^\lambda))\cap
E_{\gamma})\cup F_{\gamma}\cup D_{\lambda}.
\end{equation}
 For $\lambda$ close to $\lambda_{0}$, both the measures of $D_{\lambda}$ and
 $(\Sigma^{-}_{\lambda}\backslash B_\epsilon((z^0)^\lambda))\cap E_{\gamma}$ are close to zero. In fact, for any
$x\in(\Sigma^{-}_{\lambda}\backslash B_\epsilon((z^0)^\lambda))\cap E_{\gamma}$, we have
$$
w_{\lambda}(x)=\bar{u}_{\lambda}(x)-\bar{u}(x)
=\bar{u}_{\lambda}(x)-\bar{u}_{\lambda_{0}}(x)+
\bar{u}_{\lambda_{0}}(x)-\bar{u}(x)<0.
$$
Therefore,
$$
\bar{u}_{\lambda_{0}}(x)-\bar{u}_{\lambda}(x)>w_{\lambda_{0}}(x)>\gamma.
$$
Further,
\begin{equation}\label{n17}
((\Sigma^{-}_{\lambda}\backslash B_\epsilon((z^0)^\lambda))\cap E_{\gamma})\subset G_{\gamma}\equiv \{x\in
B_{R}|\;\;u_{\lambda_{0}}(x)-u_{\lambda}(x)>\gamma\}.
\end{equation}
By the well-known Chebyshev inequality, we have
\begin{eqnarray}\nonumber
\mu(G_{\gamma})&\leq&\frac{1}{\gamma^{p+1}}
\int_{G_{\gamma}}|u_{\lambda_{0}}(x)-u_{\lambda}(x)|^{p+1}dx\\
&\leq&
\frac{1}{\gamma^{p+1}}\int_{B_{R}}|u_{\lambda_{0}}(x)-u_{\lambda}(x)|^{p+1}dx.
\end{eqnarray}
For each fixed $\gamma$, as $\lambda$ goes to $\lambda_{0}$, the
right hand side of the above inequality goes to zero. This implies that the measure of
$(\Sigma^{-}_{\lambda}\backslash B_\epsilon((z^0)^\lambda))\cap B_{R}$ can also be made arbitrarily small.

This completes the proof of (\ref{n13.1}).

Since we can choose any direction that is perpendicular to the $x_{n}$-axis as the $x_{1}$ direction, by showing (\ref{n13.1}) or $\bar{u}$ is symmetric about the plane $T_{z_1^0}$, we have actually shown that $\bar{u}(x)$ is rotationally symmetric about the line parallel to the $x_{n}$-axis and passing through $z^{0}$ for $1< p <\frac{n+\alpha}{n-\alpha}$.
Now, for any two points $X^{1}$ and $X^{2}$ with $X^{i}=({x^{i}}^\prime,x_{n})\in R^{n-1}\times R$, $i=1,2$.  Let $z^{0}$ be the projection of $\bar{X}=\frac{X^{1}+X^{2}}{2}$ on $\partial \mathbb{R}^n_+=\{x=(x^\prime,x_{n})\in R^{n} \mid x_{n}=0\}$. Set $Y^{i}=\frac{X^{i}-z_{0}}{|X^{i}-z^{0}|^{2}}+z^{0}$, $i=1,2$. From the above arguments, one can easily see $\overline{u}(Y^{1})=\overline{u}(Y^{2})$, hence $u(X^{1})=u(X^{2}).$ This implies that $u$ is independent of $x_{1},\cdots,x_{n-1}$, which, in fact, contradicts the finiteness of the integral
$$ \int_{\mathbb{R}^n} G_{\infty} (x, y)\, y_n^2 \,u^p (y)\, dy.$$

If $u(x)=u(x_n)$ solves
\be
u(x)=\int_{\mathbb{R}^n} G_{\infty} (x, y)\, y_n^2 \,u^p (y)\, dy,
\label{n17.1}
\ee
then for each fixed $x=(x', x_n)$, let $r=|x^\prime-y^\prime|$
, $t=|x_n-y_n|$, $s=r/t$ and we have
\begin{eqnarray}\nonumber
+\infty>u(x_n)&=&
\int_{\mathbb{R}^n}\frac{A_{n,\alpha}}{|x-y|^{n-\alpha}}y_n^2 \,
u^{p}(y_n)\,  dy\\
\nonumber
&=&\int_R y_n^2 \,u^{p}(y_n)\int_{\mathbb{R}^{n-1}}
\frac{ dy'}{|x-y|^{n-\alpha}}dy_n\\
\nonumber
&=&\int_R y_n^2 \,u^{p}(y_n)\int_0^\infty
\frac{C\,r^{n-2}dr}{|r^2+t^2|^{\frac{n-\alpha}{2}}} dy_n\\
\nonumber
&=&\int_R y_n^2 \,u^{p}(y_n)\,t^{\alpha-1}\int_0^\infty
\frac{s^{n-2}\,ds}{(s^2+1)^{\frac{n-\alpha}{2}}} dy_n\\
\label{n20}
&\geq&C\int_R y_n^2 \,u^{p}(y_n)\,|x_n-y_n|^{\alpha-1}dy_n\\
\nonumber
&\geq&C\int_{R_0} y_n^2 \,u^{p}(y_n)\,(\frac{y_n}{2})^{\alpha-1}dy_n\\
\label{n18}
&=&C\int_{R_0}^\infty y_n^{1+\alpha}\,u^p(y_n)\,dy_n
\end{eqnarray}
for sufficiently large $R_0$.
The finiteness of (\ref{n18}) indicates that there exists a sequence $\{y_n^i\} \ra \infty$
as $i \ra \infty$, such that
\be\label{n19}
u^{p}(y^{i}_{n})(y^{i}_{n})^{\alpha+2}\rightarrow 0.
\ee
For any fixed $x=(0,x_{n})\in \mathbb{R}^{n-1}\times\mathbb{R}$, let $x_{n}=2R$ be large. From (\ref{n20}) we deduce that
\begin{eqnarray}
\nonumber
+\infty>u(x_{n})&\geq& C\int_1^2u^{p}(y_{n})\,y_{n}^{2}\,|x_n-y_n|^{\alpha-1}\,dy_n\\\label{n21}
&\geq&C x_{n}^{\alpha-1}.
\end{eqnarray}
Together with (\ref{n20}), (\ref{n21}) indicates that for $x_{n}=2R$ sufficiently large,
\begin{eqnarray}
\nonumber
u(x_{n})&\geq& C\int_{\frac{R}{2}}^{R}[C x_{n}^{\alpha-1}]^p \,
y_n^2\,R^{\alpha-1}\, dy_n\\
\nonumber
&=& C_{p, \alpha} R^{(\alpha-1)(p+1)+3} \\
&:=&D_{p, \alpha} x_n^{(\alpha-1)(p+1)+3}
\end{eqnarray}
Repeating the substitution process above for another $m$ times and setting $x_{n}=2R$, we arrive at
\begin{eqnarray}\label{n22}
u(x_{n})&\geq& D(m,p,\alpha)\,x_{n}^{3+3p+3p^2+\cdot\cdot\cdot+3p^m
+(p^{m+1}+p^m+\cdot\cdot\cdot+1)(\alpha-1)}\\\nonumber
&= & D(m,p,\alpha)\,x_{n}^{\frac{p^{m+1}(p-p\alpha-3)+2+\alpha}{1-p}}.
\end{eqnarray}
For $1\leq\alpha<2$, it's easy to see that as $x_{n} \ra \infty$
$$u^p(x_{n})\,x_{n}^{\alpha+2} \ra \infty.$$
This contradicts (\ref{n19}).
For $0<\alpha<1$, for sufficiently large $m$, it holds that
\begin{equation}
u^{p}(x_{n})x_{n}^{\alpha+2}\geq  D(m,p,\alpha)x_{n}^{\tau (p)}\geq D(m,p,\alpha)>0,
\end{equation}
for all $x_{n}$ sufficiently large with
\begin{eqnarray}\nonumber
&&\tau (p)\\\label{n23}
&=&3+3p+\cdot\cdot\cdot+3p^m
+(p^{m+1}+p^m+\cdot\cdot\cdot+1)(\alpha-1)+\alpha+2\\\label{n24}
&=&\frac{p^{m+2}(p-p\alpha-3)+2p+\alpha}{1-p}+2 \geq 0.
\end{eqnarray}
 Again this is a contradiction with (\ref{n19}). Hence we declare that (\ref{n17.1}) admits no positive solution.

 To prove (\ref{n24}), it suffices to show that for $m$ sufficiently large and $\alpha\in (0, 1)$,
 $$\tau^\prime(p)>0,\, \forall \,p \in (1, \frac{n+\alpha}{n-\alpha}),$$
since (\ref{n23}) shows that $\tau(1)>0$.
Indeed,
\be
\tau^\prime(p)=\frac{p^{m+1}\,\big[m\big(p(\alpha-1)+3\big)\big(p-1\big)+
6p-3\alpha \,p+2\alpha\, p^2-2p^2-6\big]
+2+\alpha}{(1-p)^2},
\ee
and $\big(p(\alpha-1)+3\big)\big(p-1\big)>0$ for $n\geq2$ and $0<\alpha<1$.
Therefore, it's true that when $m$ is large enough, $\tau^\prime(p)>0$.

This completes the proof of Theorem \ref{mthm2}.

\section{Nonexistence of positive solution for \\
$(-\lap)^{\alpha/2} u(x) = f(x_n)\, u^p(x) $}

\subsection{Equivalence between the integral equation and the
pseudo-differential equation}

Let
\be
v(x)=\int_{\mathbb{R}^n}\frac{A_{n,\alpha}}{|x-y|^{n-\alpha}}\,f(y_n) \,u^{p}(y)\,dy.
\ee
Going through exactly the same reasoning as that in section 2.1, we can show that if $u(x)$ satisfies
\be
(-\lap)^{\alpha/2} u(x) = f(x_n) \,u^{p}(x), \quad x \in \mathbb{R}^n,\label{n55}
\ee
then
$$u(x)-v(x)\equiv C\geq0.$$
We can see that $C=0$. Otherwise, for every given $x\in \mathbb{R}^n$,
\begin{eqnarray*}
\infty>u(x)&\geq & v(x)\\
&=&\int_{\mathbb{R}^n}\frac{A_{n,\alpha}}{|x-y|^{n-\alpha}}\,f(y_n)\, u^{p}(y)\,dy\\
&\geq& \int_{\mathbb{R}^n}\frac{A_{n,\alpha}}{|x-y|^{n-\alpha}}\,f(y_n) \,C^p\,dy\\
&\geq& C\int_{\mathbb{R}^n\backslash D}\frac{dy}{|x-y|^{n-\alpha}}\\
&=&\infty,
\end{eqnarray*}
with $D=\{y\in \mathbb{R}^n\mid |y_n|<1\}$.
The contradiction shows that $w=0$. Therefore,
\be
u(x)=v(x)=\int_{\mathbb{R}^n}\frac{A_{n,\alpha}}{|x-y|^{n-\alpha}}\,f(y_n) \,u^{p}(y)\,dy.\label{n54}
\ee

In the distributional sense, a solution for (\ref{n54}) is also a solution for
(\ref{n55}). Actually, for any $\phi \in C_0^\infty(\mathbb{R}^n)$,
\begin{eqnarray*}
<(-\lap)^{\frac{\alpha}{2}}u, \, \phi>&=&<\int_{\mathbb{R}^n}G_{\infty}(x,y)\,f(y_n)\,u^{p}(y)\,dy, \, (-\lap)^{\frac{\alpha}{2}}\phi(x)>\\
&=& \int_{\mathbb{R}^n}\left\{\int_{\mathbb{R}^n}G_{\infty}(x,y)\,f(y_n)\, u^{p}(y)\,dy\right\}(-\lap)^{\frac{\alpha}{2}}\phi(x)\,dx\\
&=& \int_{\mathbb{R}^n}\left\{\int_{\mathbb{R}^n}G_{\infty}(x,y) (-\lap)^{\frac{\alpha}{2}}\phi(x)dx\right\}\,f(y_n)\,u^{p}(y)\,dy\\
&=& \int_{\mathbb{R}^n}\left\{\int_{\mathbb{R}^n}\delta(x-y) \phi(x)dx\right\}\,f(y_n)\,u^{p}(y)\,dy\\
&=& \int_{\mathbb{R}^n}\phi(y)\,f(y_n)\,u^{p}(y)\,dy\\
&=& <f(y_n)\,u^{p}, \,\phi>.
\end{eqnarray*}
Hence a solution for the integral equation satisfies the
differential equation as well. This proves Theorem \ref{mthm3}.

\subsection{Nonexistence of positive solutions for the integral equation}

Performing the same Kelvin transform as defined in (\ref{n4}) on $u$ that solves
\be
u(x)=\int_{\mathbb{R}^n}G_{\infty}(x,y)\,f(y_n)\,u^{p}(y)\,dy,
\label{n35}
\ee
then $\bar{u}$, the Kelvin transform of $u$, satisfies
\begin{eqnarray}
\nonumber
\bar{u}(x)
&=&\frac{1}{|x-z^{0}|^{n-\alpha}}
\int_{\mathbb{R}^n}G_{\infty}\left(\frac{x-z^{0}}
{|x-z^{0}|^{2}}+z^{0},y \right)\,f(y_n)\, u^{p}(y)\,dy \nonumber \\
&=&\int_{\mathbb{R}^n}G_{\infty}(x,y)
\frac{1}{|x-y|^{n-\alpha}}\frac{f(\frac{y_n}{|y-z^{0}|^2}) \, \bar{u}^{p}(y)}{|y-z^{0}|^{\beta}}dy,
\label{n36}
\end{eqnarray}
for all $x\in \mathbb{R}^n\backslash B_{\epsilon}(z^{0})$ and
$\epsilon>0$, where $\beta=(n-\alpha)(\tau-p)\geq 0$ and $\tau=\frac{n+\alpha}{n-\alpha}$.
By (\ref{n35}), we have
\begin{eqnarray}
&&\bar{u}(x)-\bar{u}_{\lambda}(x)\nonumber\\
&=&\int_{\Sigma_{\lambda}}\left[G_{\infty}(x,y)-G_{\infty}(x^{\lambda},y)\right]
\frac{f(\frac{y_n}{|y-z^{0}|^2})
 \,\bar{u}^{p}(y)}{|y-z^{0}|^{\beta}}dy \nonumber\\
  &+&\int_{\Sigma_{\lambda}}
\left[G_{\infty}(x,y^{\lambda})-G_{\infty}
(x^{\lambda},y^{\lambda})\right]\frac{f(\frac{y_n}{|y^\lambda-z^{0}|^2}) \, \bar{u}^{p}_{\lambda}(y)}{|y^{\lambda}-z^{0}|^{\beta}}dy\nonumber\\
&=&\int_{\Sigma_{\lambda}}\left[
G_{\infty}(x,y)-G_{\infty}
(x^{\lambda},y)\right]\left[\frac{f(\frac{y_n}{|y-z^{0}|^2})
\,\bar{u}^{p}(y)}{|y-z^{0}|
^{\beta}}-
\frac{f(\frac{y_n}{|y^\lambda-z^{0}|^2})\,\bar{u}_{\lambda}^{p}(y)}
{|y^{\lambda}-z^{0}|^{\beta}}\right]dy.\nonumber\\
\label{n39}
\end{eqnarray}

Use the same $T_\lambda$, $
\Sigma_\lambda$ and $
x^{\lambda}$ as introduced at the beginning of section 2.1. Set
$$u_\lambda(x)=u(x^\lambda) \mbox{ and }w_\lambda(x)=\bar{u}_\lambda(x)-\bar{u}(x).$$

As before, we work on a starting point for the moving of the planes in step 1
by proving that for $\lambda$ sufficiently negative,
\be
w_\lambda(x)\geq0, \mbox{ a.e. in } \Sigma_\lambda.
\label{n37}
\ee
In step 2, one shall see that the moving plane will not stop until it arrives at the limiting position $\lambda=z_1^0$. Then by moving the planes from the very right to the left, we will have
\be
w_\lambda(x)\leq0, \mbox{ a.e. in } \Sigma_\lambda,
\label{n38}
\ee
for any $\lambda\geq z_1^0$. Together with (\ref{n37}), it yields that
\be
w_{z_1^0}\equiv0,  \mbox{ a.e. in } \Sigma_{z_1^0}.
\ee

Here comes the proof in details.

\emph{Step 1.} \emph{Start moving the planes from $-\infty$ to the right as long as (\ref{n37}) holds.}

For any $\epsilon>0$, define
$$
\Sigma_{\lambda}^{-}=\{x\in \Sigma_{\lambda}\backslash B_{\epsilon}((z^{0})^{\lambda})\mid w_{\lambda}(x)<0\}.
$$
We show that for $\lambda$ sufficiently negative,
$\Sigma_{\lambda}^{-}$ must be measure zero.

It follows from the Mean Value Theorem and (\ref{n39}) that for  $x\in\Sigma_{\lambda}^{-}$ and $\lambda$ sufficiently negative, we have
\begin{eqnarray}
\nonumber
0&<&\bar{u}(x)-\bar{u}_{\lambda}(x)\\
\nonumber
&= &\int_{\Sigma_{\lambda}}\left[
G_{\infty}(x,y)-G_{\infty}
(x^{\lambda},y)\right]\left[\frac{f(\frac{y_n}{|y-z^{0}|^2})
\,\bar{u}^{p}(y)}{|y-z^{0}|
^{\beta}}-
\frac{f(\frac{y_n}{|y^\lambda-z^{0}|^2})\,\bar{u}_{\lambda}^{p}(y)}
{|y^{\lambda}-z^{0}|^{\beta}}\right]dy\\
\nonumber
&=&\int_{\Sigma_{\lambda^-}}\left[
G_{\infty}(x,y)-G_{\infty}
(x^{\lambda},y)\right]\left[\frac{f(\frac{y_n}{|y-z^{0}|^2})
\,\bar{u}^{p}(y)}{|y-z^{0}|
^{\beta}}-
\frac{f(\frac{y_n}{|y^\lambda-z^{0}|^2})\,\bar{u}_{\lambda}^{p}(y)}
{|y^{\lambda}-z^{0}|^{\beta}}\right]dy\\
\nonumber
&\;\;+&\int_{\Sigma_{\lambda}\backslash \Sigma_{\lambda^-}}\left[
G_{\infty}(x,y)-G_{\infty}
(x^{\lambda},y)\right]\left[\frac{f(\frac{y_n}{|y-z^{0}|^2})
\,\bar{u}^{p}(y)}{|y-z^{0}|
^{\beta}}-
\frac{f(\frac{y_n}{|y^\lambda-z^{0}|^2})\,\bar{u}_{\lambda}^{p}(y)}
{|y^{\lambda}-z^{0}|^{\beta}}\right]dy\\
\nonumber
&\leq& \int_{\Sigma_{\lambda^-}}\left[
G_{\infty}(x,y)-G_{\infty}
(x^{\lambda},y)\right]\left[\frac{f(\frac{y_n}{|y-z^{0}|^2})
\,\bar{u}^{p}(y)}{|y-z^{0}|
^{\beta}}-
\frac{f(\frac{y_n}{|y^\lambda-z^{0}|^2})\,\bar{u}_{\lambda}^{p}(y)}
{|y^{\lambda}-z^{0}|^{\beta}}\right]dy\\
\nonumber
&\leq&\int_{\Sigma^{-}_{\lambda}}\left[
G_{\infty}(x,y)-G_{\infty}(x^{\lambda},y)\right]
f(\frac{y_n}{|y-z^{0}|^2})\,
\frac{\bar{u}^{p}(y)-\bar{u}_{\lambda}^{p}(y)}
{|y-z^{0}|^{\beta}}dy\\
\nonumber
&\leq&p\int_{\Sigma^{-}_{\lambda}}
G_{\infty}(x,y)\,f(\frac{y_n}{|y-z^{0}|^2})\,
\frac{\bar{u}^{p-1}(y)}{|y-z^{0}|^{\beta}}
[\bar{u}(y)-\bar{u}_{\lambda}(y)]dy\\
\label{n40}
&=&\int_{\Sigma^{-}_{\lambda}}
\frac{C}{|x-y|^{n-\alpha}}\,f(\frac{y_n}{|y-z^{0}|^2})\,
\frac{\bar{u}^{p-1}(y)}{|y-z^{0}|^{\beta}}
(\bar{u}(y)-\bar{u}_{\lambda}(y))dy.
\end{eqnarray}

  Let $\Omega$ be any domain that is a positive distance away from $z^0$. Since $u$ is locally bounded, we have
\begin{equation}\label{n41}
\int_{\Omega}\left[\frac{f(\frac{y_n}{|y^\lambda-z^0|^2})
\bar{u}^{p-1}(y)}{|y-z^{0}|^{\beta}}\right]
^{\frac{n}{\alpha}}dy<
C\,\int_{\Omega}\left[\frac{\bar{u}^{p-1}(y)}{|y-z^{0}|^{\beta}}\right]
^{\frac{n}{\alpha}}dy<\infty.
\end{equation}

For any $q>\frac{n}{n-\alpha}$, applying the Hardy-Littlewood-Sobolev inequality (\ref{n8}) and H\"{o}lder
inequality to (\ref{n40}) yields
\begin{eqnarray}\nonumber
\|w_{\lambda}\|_{L^{q}(\Sigma_{\lambda}^{-})}&\leq&
C\|f(\frac{y_n}{|y^\lambda-z^0|^2})\,
\frac{\bar{u}^{p-1}}{|y-z^{0}|^{\beta}}w_{\lambda}\|_{L^{\frac{nq}{n+\alpha
q}}(\Sigma_{\lambda}^{-})}\\\nonumber
&\leq& C\|f(\frac{y_n}{|y^\lambda-z^0|^2})\,
\frac{\bar{u}^{p-1}}{|y-z^{0}|^{\beta}}\|_{L^{\frac{n}{\alpha}}
(\Sigma_{\lambda}^{-})}
\|w_{\lambda}\|_{L^{q}(\Sigma_{\lambda}^{-})}\\
\label{n42}
&\leq& C\|\frac{\bar{u}^{p-1}}{|y-z^{0}|^{\beta}}\|_{L^{\frac{n}{\alpha}}
(\Sigma_{\lambda}^{-})}
\|w_{\lambda}\|_{L^{q}(\Sigma_{\lambda}^{-})}.
\end{eqnarray}
When $N$ is sufficiently large, (\ref{n41}) indicates that for $\lambda\leq-N$,
\begin{equation}\label{n43}
C\left
\{\int_{\Sigma_{\lambda}^{-}}
[\frac{\bar{u}^{p-1}}{|y-z^{0}|^{\beta}}]^{\frac{n}{\alpha}}dy\right\}
^{\frac{\alpha}{n}}\leq\frac{1}{2}.
\end{equation}
A combination of (\ref{n42}) and (\ref{n43}) gives
$$
\|w_{\lambda}\|_{L^{q}(\Sigma_{\lambda}^{-})}=0.
$$
Therefore $\Sigma_{\lambda}^{-}$ must be measure zero, that is,
\begin{equation}\label{n44}
w_{\lambda}(x)\geq 0, \;\;a.e.\;x\in \Sigma_{\lambda}.
\end{equation}

\emph{Step 2.} {\em Move the plane to the limiting position to
derive symmetry}.

Let
$$
\lambda_{0}=\sup\{\lambda\leq z_{1}^{0}\mid w_{\rho}(x)\geq 0,\; \rho\leq\lambda,\;
\forall x\in\Sigma_{\rho}\}.
$$
We show that
$$\lambda_{0}=z_{1}^{0}$$ and
$\bar{u}(x)$ is symmetric about the plane
$T_{\lambda_{0}}$, i.e.
\begin{equation}\label{n45}
w_{\lambda_{0}}\equiv 0,\;\;\mbox{ a.e. in }
\Sigma_{\lambda_{0}}.
\end{equation}

Despite that
$\beta$ takes a different value here, the rest of the proof will be the same as that for $(-\lap)^{\alpha/2}u(x)=x_n^2\,u^p(x)$, because $f(\frac{y_n}{|y-z^0|^2})
\bar{u}^{p-1}(y)$ can be controlled by $\bar{u}^{p-1}(y)$ for
$0<|\lambda-z_{1}^{0}|<\epsilon$ with fixed $\epsilon$.

The arbitrariness of the choice for the $x_1$ direction contributes to
$\bar{u}(x)$'s rotational symmetry about the line parallel to the $x_{n}$-axis and passing through $z^{0}$. This implies that $u$ relies on its $x_n$th
variable only. However, if $u(x)=u(x_n)$ solves
\be
u(x)=\int_{\mathbb{R}^n} G_{\infty} (x, y)\, f(y_n) \,u^p (y)\, dy,
\label{n46}
\ee
then for each fixed $x=(x', x_n)$, let $r=|x^\prime-y^\prime|$
, $t=|x_n-y_n|$, $s=r/t$ and it gives
\begin{eqnarray}\nonumber
+\infty>u(x_n)&=&
\int_{\mathbb{R}^n}\frac{A_{n,\alpha}}{|x-y|^{n-\alpha}}f(y_n) \,
u^{p}(y_n)\,  dy\\
\nonumber
&=&\int_R f(y_n) \,u^{p}(y_n)\int_{\mathbb{R}^{n-1}}
\frac{ dy'}{|x-y|^{n-\alpha}}dy_n\\
\nonumber
&=&\int_R f(y_n) \,u^{p}(y_n)\int_0^\infty
\frac{C\,r^{n-2}dr}{|r^2+t^2|^{\frac{n-\alpha}{2}}} dy_n\\
\nonumber
&=&\int_R f(y_n) \,u^{p}(y_n)\,t^{\alpha-1}\int_0^\infty
\frac{s^{n-2}\,ds}{(s^2+1)^{\frac{n-\alpha}{2}}} dy_n\\
\label{n47}
&\geq&C\int_R f(y_n) \,u^{p}(y_n)\,|x_n-y_n|^{\alpha-1}dy_n\\
\nonumber
&\geq&C\int_{R_0} f(y_n) \,u^{p}(y_n)\,(\frac{y_n}{2})^{\alpha-1}dy_n\\
\label{n48}
&=&C\int_{R_0}^\infty f(y_n)\,y_n^{\alpha-1}\,u^p(y_n)\,dy_n
\end{eqnarray}
for sufficiently large $R_0$.
The finiteness of (\ref{n48}) indicates that there exists a sequence $\{y_n^i\} \ra \infty$
as $i \ra \infty$, such that
\be
f(y_n^i)\,u^{p}(y^{i}_{n})(y^{i}_{n})^{\alpha}\rightarrow 0.
\ee
Due to the monotonicity of $f(y_n)$, we further deduce that
\be\label{n49}
u^{p}(y^{i}_{n})(y^{i}_{n})^{\alpha}\rightarrow 0.
\ee

For any fixed $x=(0,x_{n})\in \mathbb{R}^{n-1}\times\mathbb{R}$, let $x_{n}=2R$ be large. From (\ref{n47}) we deduce that
\begin{eqnarray}
\nonumber
+\infty>u(x_{n})&\geq& C\int_1^2u^{p}(y_{n})\,f(y_{n})\,|x_n-y_n|^{\alpha-1}\,dy_n\\\label{n50}
&\geq&C x_{n}^{\alpha-1}.
\end{eqnarray}
Together with (\ref{n47}), (\ref{n50}) indicates that for $x_{n}=2R$ sufficiently large,
\begin{eqnarray}
\nonumber
u(x_{n})&\geq& C\int_{\frac{R}{2}}^{R}[C x_{n}^{\alpha-1}]^p \,
f(y_{n})\,R^{\alpha-1}\, dy_n\\
\nonumber
&=& C_{p, \alpha} R^{(\alpha-1)(p+1)+1} \\
&:=&D_{p, \alpha} x_n^{\alpha+\alpha \,p-p}
\end{eqnarray}
Going through the substitution process above for another $m$ times and setting $x_{n}=2R$, it gives
\be\label{n51}
u(x_{n})\geq D(m,p,\alpha)\,x_{n}^{\alpha+\alpha\, p+\alpha\,
 p^2+\cdot\cdot\cdot+\alpha\, p^{m+1}-p^{m+1}}.
\ee
For $1\leq\alpha<2$, it's easy to see that as $x_{n} \ra \infty$
$$u^p(x_{n})\,x_{n}^{\alpha} \ra \infty.$$
This contradicts (\ref{n49}).
For $0<\alpha<1$, for sufficiently large $m$, it holds that
\begin{equation}
u^{p}(x_{n})x_{n}^{\alpha}\geq  D(m,p,\alpha)x_{n}^{\tau (p)}\geq D(m,p,\alpha)>0,
\end{equation}
for all $x_{n}$ sufficiently large with
\begin{eqnarray}\nonumber
&&\tau (p)\\\label{n52}
&=&\alpha+\alpha\, p+\alpha\,
 p^2+\cdot\cdot\cdot+\alpha\, p^{m+1}-p^{m+1}+\alpha\\\label{n53}
&=&\frac{\alpha-p^{m+2}-(\alpha-1)\,p^{m+3}}{1-p} \geq 0.
\end{eqnarray}
 Again this is a contradiction with (\ref{n49}). Hence we declare that (\ref{n35}) admits no positive solution.
 To prove (\ref{n53}), it suffices to show that for $m$ sufficiently large and $\alpha\in (0, 1)$,
 $$\tau^\prime(p)>0,\, \forall \,p \in (1, \frac{n+\alpha}{n-\alpha}),$$
since (\ref{n52}) shows that $\tau(1)>0$.
Indeed,
\be
\tau^\prime(p)=\frac{p^{m+1}\,\big[m\big(p(\alpha-1)+1\big)\big(p-1\big)+
4p^2(\alpha-1)+3p(2-\alpha)-2\big]
-\alpha}{(1-p)^2},
\ee
and $\big(p(\alpha-1)+1\big)\big(p-1\big)>0$ for $n\geq2$ and $0<\alpha<1$.
Therefore, it's true that when $m$ is large enough, $\tau^\prime(p)>0$.

This completes the proof of Theorem \ref{mthm4}.

{\em the Author's Address and E-mail:}
\medskip

Yan Li

Department of Mathematical Sciences

Yeshiva University

New York, NY, 10033 USA

yali3@mail.yu.edu


\begin{thebibliography}{CL}

\bibitem[CC]{CC} L. Cao and W. Chen, \,\,Liouville type theorems for poly-harmonic Navier problems, \,\,
    Disc. Cont. Dyn. Sys. 33 (2013) 3937-3955.

 \bibitem[CDL]{CDL}  W. Chen, L. D' Ambrosio and Y. Li, \,\, Some Liouville theorems for
the fractional Laplacian, \,\,Nonlinear Anal: Theory, Methods \& Appl. in
press, doi:10.1016/j.na.2014.11.003.


\bibitem[CFY]{CFY}  W. Chen, Y. Fang  and R. Yang, \,\, Semilinear equations involving the fractional Laplacian on domains, \,\,
     Advances in Math, in press.

\bibitem[CL1]{CL1} W. Chen and C. Li, \,\, Regularity of solutions for a system of integral equation,
 \,\, Comm. Pure Appl. Anal, 4(2005) 1-8.

\bibitem[CL2]{CL2} W. Chen and C. Li, \,\, Methods on Nonlinear Elliptic Equations, \,\, AIMS. Ser. Differ. Equ. Dyn.
Syst, vol.4 2010.

\bibitem[CL3]{CL3} W. Chen and C. Li,  \,\,Radial symmetry of solutions for some integral
systems of Wolff type, \,\, Disc. Cont. Dyn. Sys. 30(2011) 1083-1093.

\bibitem[CLO]{CLO} W. Chen, C. Li, and B. Ou, \,\, Qualitative properties of solutions for
an integral equation, \,\,Disc. Cont. Dyn. Sys. 12(2005) 347-354.

\bibitem[CS]{CS} L. Caffarelli and L. Silvestre, \,\, An extension problem related to the fractional Laplacian, \,\,
Comm. PDE. 32(2007) 1245-1260.

\bibitem[Ku]{Ku} T. Kulczycki, \,\, Properties of Green function of symmetric stable processes,
    \,\, Prob. Math. Stat. 17(1997) 339-364.


\bibitem[LZ]{LZ} D. Li and R. Zhuo, \,\,   An integral equation on half space, \,\, Proc. AMS. 138(2010) 2779-2791.

\bibitem[MC]{MC} L. Ma and D. Chen, \,\,A Liouville type theorem for an integral system, \,\, Comm. Pure Appl. Anal. 5(2006) 855-859.

\bibitem[RFB]{RFB} R. Zhuo, F. Li and B. Lv, \,\,Liouville type theorems for Schrödinger system with Navier boundary conditions in a half space,\,\,
     Comm. Pure Appl. Anal. 13(2014) 977-990.

\bibitem[Si]{Si} L. Silvestre, \,\, Regularity of the obstacle problem for
a fractional power of the Laplace operator, \,\, Comm. Pure Appl. Math. 60(2007) 67-112.

\bibitem[ZCCY]{ZCCY} R. Zhuo, W. Chen, X. Cui and Z. Yuan, \,\, A Liouville theorem for the fractional Laplacian, \,\,
     Accepted by Disc. Cont. Dyn. Sys.

\bibitem[ZW]{ZW} Yonggang  Zhao and Mingxin  Wang,\,\,
An integral equation involving Bessel potentials on half space, \,\,
Comm. Pure Appl. Anal. 14(2015) 1534-0392.


\end{thebibliography}
\end{document}